\documentclass[english]{article}

\usepackage{mathrsfs}%
\usepackage[T1]{fontenc}
\usepackage{marvosym}
\usepackage{babel}
\usepackage{amsmath,amsfonts,amsthm}%
\usepackage{amsmath}%
\usepackage{amsfonts}%
\usepackage{amssymb}%
\usepackage{graphicx}
\newcommand{\lp}{\mathbb{L}^{p}}

\newcommand{\lpc}{\mathbb{L}^{p}_c}
\newcommand{\linf}{\mathbb{L}^{\infty}}

\newcommand{\N}{\mathbb{N}}

\newcommand{\E}{\mathbb{E}}

\newcommand{\PP}{{\cal P}}

\newtheorem{theorem}{Theorem}

\newtheorem{corollary}[theorem]{Corollary}

\newtheorem{definition}[theorem]{Definition}
\newtheorem{example}[theorem]{Example}

\newtheorem{lemma}[theorem]{Lemma}

\newtheorem{proposition}[theorem]{Proposition}
\newtheorem{remark}[theorem]{Remark}

\newcommand{\vep}{\varepsilon}
\newcommand{\tX}{\tilde{X}}
\begin{document}

\title{Function spaces and capacity related to a Sublinear
Expectation: application to G-Brownian Motion Paths}

\author{Laurent Denis $\cdot$ Mingshang Hu $\cdot$ Shige
Peng} \maketitle \makeatletter
\newskip\@footindent
\@footindent=0em
\renewcommand\footnoterule{\kern-3\p@ \hrule width 0.4\columnwidth \kern 2.6\p@}

\long\def\@makefntext#1{\@setpar{\@@par\@tempdima \hsize
\advance\@tempdima-\@footindent
\parshape \@ne \@footindent \@tempdima}\par
\noindent \hbox to \z@{\hss\@thefnmark\hspace{0.2em}}#1}
\renewcommand\thefootnote{\myfootnotestyle{\arabic{footnote}}}
\def\@makefnmark{\hbox{\@thefnmark}}
\newcommand\myfootnotestyle[1]{\ifcase#1 \or £¿ \or £¿ \or £¿ \else *\fi\relax}
\makeatother

\renewcommand\thefootnote{}\footnote{L. Denis
\newline  D\'{e}partement de Math\'{e}matiques, Equipe
``Analyse et Probabilit\'{e}s", \newline Universit\'{e}
d'Evry-Val-d'Essonne,\newline Boulevard F. Mitterrand, 91 025 EVRY
Cedex, France \newline e-mail: ldenis@univ-evry.fr
\newline
\newline M. Hu (\Letter)  \newline School of Mathematics, Shandong
University,\newline Jinan, 250100, China \newline e-mail:
humingshang@mail.sdu.edu.cn
\newline
\newline S. Peng \newline
School of Mathematics, Shandong University,\newline Jinan, 250100,
China \newline School of Mathematics, Fudan University,\newline
Shanghai, 200433, China \newline e-mail: peng@sdu.edu.cn }
%"Title of the paper"

%\begin{abstract}
%\begin{keyword}[class=AMS]
%\kwd[Primary ]{60H05, 60G44}  \kwd[; secondary ]{31C15}
%\end{keyword}
%\begin{keyword}
%\kwd{}
%%\kwd{}
%\end{keyword}

\noindent\textbf{Abstract} In this paper we give some basic and
important properties of several typical Banach spaces of functions
of $G$-Brownian motion paths induced by a sublinear
expectation--$G$-expectation. Many results can be also applied to
more general situations.  A generalized version of Kolmogorov's
criterion for continuous modification of a stochastic process is
also obtained. The results can be applied in continuous time dynamic
and coherent risk measures in finance in particular for
path-dependence risky positions under situations of volatility model
uncertainty.

\  \  \


\section{Introduction}

%A Wiener probability measure $P_W$ is defined on the space of
%continuous paths
%$(\Omega,\mathscr{F})=(C[0,\infty),\mathscr{B}(C[0,\infty)))$ in
%which the canonical process $B_t(\omega)=\omega_t$, $t\geq 0$,
%$\omega\in\Omega$ becomes a Brownian motion. $B$ is a continuous
%process with zero-mean, independent and stationary increments, which
%can be proved to be normally distributed. The expectation $
%E_{P_W}[\cdot]$ defined on $L^1_{P_W}(\Omega)$, the space of all
%$P_W$-integrable $\mathscr{F}$-measurable random variables, forms a
%constant-preserving, monotone and linear functional.
%
%Since the paper \cite{ADEH} on coherent risk measures, people are
%more and more interested in  sublinear expectations (or more
%generally, convex expectations, see \cite{FoSch,FrRo}).
%
%A very interesting phenomenon is that there exists a sublinear
%expectation $\mathbb{E}[\cdot]$ on some well-defined space
%$L^1_G(\Omega)$, which plays the role of $L^1_{P_W}(\Omega)$,
% for dealing with uncertainty
%volatility model (UVM) initialed in \cite{avel,terry} (see also
%\cite{atlan,elkqu}).

How to measure the risk of financial losses in a financial market is
still a  challenging  problem. In a seminal paper \cite{ADEH} a
basic notion of coherent risk measures was introduced: Let
$\mathcal{H}$ be a linear space of financial losses, considered as a
space of random variables. A coherent risk measure
$\mathbb{E}:\mathcal{H\mapsto}\mathbb{R}$ is a real valued (monetary
value) functional with the properties of constant preserving (called
cash invariance), monotonicity, convexity and positive homogeneity.
Namely, a coherent risk measure is in fact a sublinear expectation
$\mathbb{E}$ defined on $\mathcal{H}$ (see Theorem \ref{THM2} or
Definition \ref{def3.1}). It was proved that a sublinear expectation
has the following representation (see \cite{ADEH, Delb, FoSch}):
There exists a family of linear expectations $\left \{
E_{\theta}\right \} _{\theta \in \Theta}$ such that
\[
\mathbb{E}[X]=\max_{\theta \in \Theta}E_{\theta}[X],\  \ X\in
\mathcal{H}.
\]
The meaning in economics of this representation is that the risk
measure $\mathbb{E}$ is in fact the robust super-expectation over
the family of uncertainty of linear expectations $\left \{
E_{\theta}\right \} _{\theta \in \Theta}$.

As an example, let us consider a typical situation in a financial
market where the price of a stock satisfies the equation
$dS_{t}=S_{t}(\gamma_{t}^{\theta }dt+\sigma_{t}^{\theta}dW_{t})$,
where $W$ is a standard Brownian motion, and
$(\gamma_{t}^{\theta},\sigma_{t}^{\theta})_{t\geq0}$, $\theta \in
\Theta$, are unknown processes parameterized by $\theta \in \Theta$.
A financial loss $X$ is formulated as a given random variable
depending on the path of $S$ or, equivalently, on the path of
$B_{t}^{\theta}=\int_{0}^{t}(\gamma_{s}^{\theta
}ds+\sigma_{s}^{\theta}dW_{s})$. For each fixed  $\theta \in
\Theta$, let $P_{\theta}$ be the probability measure on the space of
continuous paths
$(\Omega,\mathcal{F})=(C(0,\infty),\mathcal{B}(C(0,\infty)))$
induced by $(B_{t}^{\theta})_{t\geq 0}$, and $E_{\theta}$ the
corresponding expectation. The risk measure of $X$
under the above uncertainty is formulated as%
\[
\mathbb{E}[X]=\sup_{\theta \in \Theta}E_{\theta}[X].
\]
A typical situation is that each probability $P_{\theta}$ is
absolutely continuous with respect to the `reference measure'
$P_{0}$ corresponding to the case when $B_t^{\theta_0}\equiv W_t$.
In this case the uncertainty can only comes from $\gamma^{\theta}$
and thus is called drift uncertainty. Chen and Epstein \cite{CE}
proposed to use $g$-expectation, (small $g$, introduced in
\cite{Peng1997}) for a robust valuation:
$\mathbb{E}_{g}[X]=\sup_{\theta \in \Theta}E_{\theta}[X]$. It was
also proved (see \cite{EPQ, CE}) that this corresponds to the case
where the uncertain drift has the form $\{ \gamma_{s}^{\theta}\in
K,s\geq0\}$, for some $K\subset {\mathbb R}$. \cite{Roazza2003,
F-RG1} proposed to use $g$-expectation as a time consistent risk
measure.  \cite{DPR} proved that any coherent and time consistent
risk measure absolutely continuous with respect to $P_{0}$ can be
approximated by a $g$-expectation.

But in finance there is an important situation called `volatility
uncertainty' in which the uncertainty comes from the
\textquotedblleft volatility coefficient\textquotedblright \
$\{\sigma^{\theta}, \theta \in \Theta\}
$. A major difficulty here is that the probabilities $\{P_{\theta}%
\}_{\theta \in \Theta}$ are mutually singular and thus the
corresponding $\mathbb{E}$ cannot be dominated by any
$g$-expectation. This type of uncertainty was initially studied by
Avellaneda, Levy and Paras \cite{avel} and Lyons \cite{terry}, for
the superhedging of European options with payoffs depending only on
the terminal value $B_{T}^{\theta}$, the discrete-time case has been
also studied in [10]. But for the superhedging of a general
path-dependence option, the difficulty was dramatically increased.
This situation was studied independently  by \cite{Peng2004} and
\cite{Denmar} with very different approaches. Motivated by the
problem of coherent risk measures under the volatility uncertainty,
\cite{Peng2006a,Peng2006b} introduced a sublinear expectation on a
well-defined space $L_{G}^{1}(\Omega)$ under which the increments of
the canonical process $(B_t)_{t\geq 0}$ are zero-mean, independent
and stationary and can be proved to be `$G$-normally
distributed'(see \cite{Peng2007b}). This type of processes is called
`$G$-Brownian motion' and the corresponding sublinear expectation
$\mathbb{E}[\cdot]$ is called `$G$-expectation' (capital $G$).
Recently, we have discovered a strong link between the framework of
\cite{Peng2004,Peng2005,Peng2006a, Peng2006b} and the one introduced
in \cite{Denmar}.

A well-known and fundamentally important fact in probability theory
is that the linear space $L^1_{P_0}(\Omega)$ coincides with the
$E_{P_0}[|\cdot|]$-norm completion of the space of bounded and
continuous functions $C_b(\Omega)$ or bounded and
$\mathscr{F}$-measurable functions $B_b(\Omega)$, or even smaller
one, the space $L_{ip}(\Omega)\subset C_b(\Omega)$ of bounded and
Lipschitz cylinder functions (see Section 3 for its definition).

Similar problems arise in the theory of $G$-Brownian motion: The
space $L^1_G(\Omega)$ is defined as the $\mathbb{E}[|\cdot|]$-norm
completion of  $L_{ip}(\Omega)$. Can we prove that each element
$X\in L^1_G(\Omega)$ can be identified as an element of
$\mathbb{L}^1$, the space of all $\mathscr{F}$-measurable random
variables $X$ such that $\mathbb{E}[|X|]<\infty$? Furthermore, what
is the relation between the $\mathbb{E}[|\cdot|]$-norm completions
of $B_b(\Omega)$, $C_b(\Omega)$ and $L_{ip}(\Omega)$?

In this paper we give an affirmative answer to the first problem.
For the second problem, we will prove that, in fact, the
$\mathbb{E}[|\cdot|]$-completions of $L_{ip}(\Omega)$ and
$C_b(\Omega)$ are the same, but they are strict subspace of the
$\mathbb{E}[|\cdot|]$-completion of $B_b(\Omega)$. Moreover the
later one is, again, a strict subspace of $\mathbb{L}^1$.

In this paper a weakly compact family $\mathscr{P}$ of probability
measures on $(\Omega,\mathscr{F})$ is constructed so that the $G$
expectation is the upper expectation of $\mathscr{P}$, i.e.:
$$
\mathbb{E}[X]=\sup_{P\in \mathscr{P}}E_P[X]\, \, \,  \textrm{ for
each }\, X\in L_{ip}(\Omega).
$$
Following \cite{HuSt}, we define the corresponding regular Choquet
capacity:
$$
c(A):=\sup_{P\in \mathscr{P}}P(A),\, \, \, A\in\mathscr{B}(\Omega).
$$
We then can prove that each element $X\in L^1_G(\Omega)$ has a
$c$-quasi continuous version on $\Omega$. Moreover we have
$C_b(\Omega)\subset L^1_G(\Omega)\subset\mathbb{L}^1$ (see also
\cite{fedlp,Denmar} for a different approach).

This paper is organized as follows: in Section 2, we use a family of
probability measures $\mathscr{P}$ on $(\Omega,\mathscr{B}(\Omega))$
to define a sublinear expectation as the upper expectation of
$\mathscr{P}$, as well as the related capacity, especially, we use a
weakly compact family of probability measures to define the
corresponding regular sublinear expectation and regular capacity.
Here $\Omega$ is assumed to be a general complete separable metric
space of which $C(0,\infty)$ and $D(0,\infty)$ (the path space of
c\`{a}dl\`{a}g processes) are typical examples. Each element of
$\mathbb{E}[|\cdot|^p]^{1/p}$-completion of $C_b(\Omega)$ is proved
to has a quasi-continuous version. Concrete characterizations of
completions of different function spaces are given. As a by-product,
we  obtain a generalized  version of Kolmogorov's criterion for
continuous modification of a stochastic process. In Section 3, we
let $\Omega=C^d[0,\infty)$ and use a method of stochastic control to
prove that $G$-expectation is an upper expectation associated to a
weakly compact family $\mathscr{P}$ and then apply the results of
Section 2 to the $G$-expectation and the corresponding functional
spaces.

\section{Integration theory associated to an upper probability}

Let $\Omega$ be a complete separable metric space equipped with the
distance $d$, $\mathscr{B}(\Omega)$ the Borel $\sigma$-algebra of
$\Omega$ and $\mathcal{M}$ the collection of all probability
measures on $(\Omega ,\mathscr{B}(\Omega))$.

\begin{itemize}
\item $L^{0}(\Omega)$: the space of all $\mathscr{B}(\Omega)$-measurable real functions;

\item $B_{b}(\Omega)$: all bounded functions in $L^{0}(\Omega)$;

\item $C_{b}(\Omega)$: all continuous functions in $B_{b}(\Omega)$.
\end{itemize}

All along this section, we consider a given subset
$\mathscr{P}\subseteq \mathcal{M}$.

\subsection{Capacity associated to $\mathscr{P}$\  \ }

We denote
\[
c(A):=\sup_{P\in \mathscr{P}} P(A),\  \  \  \ A\in
\mathscr{B}(\Omega).
\]

One can easily verify the following theorem.

\begin{theorem}
The set function $c(\cdot)$ is a Choquet capacity, i.e. (see
\cite{Cho,del}),
\begin{enumerate}
\item $0\leq c(A)\leq1,\  \  \forall A\subset \Omega$.

\item If $A\subset B$, then $c(A)\leq c(B)$.

\item If $(A_{n})_{n=1}^{\infty}$ is a
sequence in $\mathscr{B}(\Omega)$, then $c(\cup A_{n})\leq \sum c(A_{n}%
)$.

\item If $(A_{n})_{n=1}^{\infty}$ is an increasing
sequence in $\mathscr{B}(\Omega)$: $A_{n}\uparrow A=\cup A_{n}$,
then $c(\cup A_{n})=\lim _{n\rightarrow \infty}c(A_{n})$.
\end{enumerate}
\end{theorem}

Furthermore, we have

\begin{theorem}
For each $A\in \mathscr{B}(\Omega)$, we have
\[
c(A) =\sup \{ c(K):\ K\makebox{ compact }\ K\subset A\}.
\]
\end{theorem}

\begin{proof} It is simply because
\[
c(A) =\sup_{P\in \mathscr{P}}\sup_{\substack{K\,
\text{compact}\\K\subset A}}P(K)=\sup_{\substack{K\,
\text{compact}\\K\subset A}}\sup_{P\in
\mathscr{P}}P(K)=\sup_{\substack{K\, \text{compact}\\K\subset
A}}c(K).
\]
\end{proof}

\begin{definition}
We use the standard capacity-related vocabulary: a set $A$ is
\textit{polar} if $c(A)=0$ and a property holds ``quasi-surely''
(q.s.) if it holds outside a polar set.
\end{definition}

\begin{remark}
In other words, $A\in \mathscr{B} (\Omega)$ is polar if and only if
$P(A)=0$ for any $P\in \mathscr{P}$.
\end{remark}

We also have in a trivial way a Borel-Cantelli Lemma.

\begin{lemma}
{\label{BorelC}} Let $(A_{n} )_{n\in \mathbb{N}}$ be a sequence of
Borel sets such that
\[
\sum_{n=1}^{\infty} c(A_{n} )<\infty.
\]
Then $\limsup_{n\rightarrow\infty}A_{n}$ is polar .
\end{lemma}

\begin{proof} Applying the Borel-Cantelli Lemma under each probability
$P\in \mathscr{P}$.
\end{proof}

The following theorem is Prohorov's theorem.

\begin{theorem} \label{newth6}
$\mathscr{P}$ is relatively compact if and only if for each
$\varepsilon >0$, there exists a compact set $K$ such that
$c(K^{c})<\varepsilon$.
\end{theorem}

The following two lemmas can be found in \cite{HuSt}.

\begin{lemma}
\label{Lemma1}$\mathscr{P}$ is relatively compact if and only if for
each sequence of closed sets $F_{n}\downarrow \emptyset$, we have
$c(F_{n})\downarrow0$.
\end{lemma}

\begin{proof} We outline the
proof for the convenience of readers.\\
\textquotedblleft$\Longrightarrow$\textquotedblright \ part: It
follows from Theorem \ref{newth6} that for each fixed
$\varepsilon>0$, there exists a compact set $K$ such that
$c(K^{c})<\varepsilon$. Note that $F_{n}\cap K \downarrow
\emptyset$, then there exists an $N>0$ such that $F_{n}\cap K
=\emptyset$ for $n \geq N$, which implies
$\lim_{n}c(F_n)<\varepsilon$. Since $\varepsilon$ can be arbitrarily
small, we obtain $c(F_{n})\downarrow0$.
\newline \textquotedblleft$\Longleftarrow$\textquotedblright \ part:
For each $\varepsilon>0$, let $(A_{i}^{k})_{i=1}^{\infty}$ be a
sequence of open balls of radius $1/k$ covering $\Omega$. Observe
that $(\cup_{i=1}^{n}A_{i}^{k})^{c}\downarrow \emptyset$, then there
exists an $n_k$ such that
$c((\cup_{i=1}^{n_k}A_{i}^{k})^{c})<\varepsilon 2^{-k}$. Set
$K=\overline{\cap_{k=1}^{\infty}\cup_{i=1}^{n_k}A_{i}^{k}}$. It is
easy to check that $K$ is compact and $c(K^{c})<\varepsilon$. Thus
by Theorem \ref{newth6} $\mathscr{P}$ is relatively compact.
\end{proof}

\begin{lemma}\label{newle8}
Let $\mathscr{P}$ be weakly compact. Then for each sequence of
closed sets $F_{n}\downarrow F$, we have $c(F_{n})\downarrow c(F)$.
\end{lemma}

\begin{proof} We outline the
proof for the convenience of readers. For each fixed
$\varepsilon>0$, by the definition of $c(F_n)$, there exists a $P_n
\in \mathscr{P}$ such that $P_n(F_n) \geq c(F_n)-\varepsilon$. Since
$\mathscr{P}$ is weakly compact, there exist $P_{n_k}$ and $P \in
\mathscr{P}$ such that $P_{n_k}$ converge weakly to $P$. Thus
\[
P(F_m) \geq \limsup_{k\rightarrow\infty}P_{n_k}(F_m) \geq
\limsup_{k\rightarrow\infty}P_{n_k}(F_{n_k}) \geq
\lim_{n\rightarrow\infty}c(F_n)-\varepsilon.
\]
Letting $m\rightarrow\infty$, we get $P(F)\geq
\lim_{n\rightarrow\infty}c(F_n)-\varepsilon$, which yields
$c(F_{n})\downarrow c(F)$.
\end{proof}

Following \cite{HuSt} (see also \cite{Delbaen,FoSch}) the upper
expectation of $\mathscr{P}$ is defined as follows: for each $X\in
L^{0}(\Omega)$ such that $E_{P}[X]$ exists for each $P\in
\mathscr{P}$,
\[
\mathbb{E} [X]=\mathbb{E} ^{\mathscr{P}}[X]:=\sup_{P\in \mathscr{P}}%
E_{P}[X].
\]

It is easy to verify

\begin{theorem}\label{THM2}
The upper expectation $\mathbb{E} [\cdot]$ of the family
$\mathscr{P}$  is a sublinear expectation on $B_{b}(\Omega)$ as well
as on $C_{b}(\Omega)$, i.e.,
\begin{enumerate}
\item for all $X,Y$ in $B_{b}(\Omega)$, $X\geq
Y\Longrightarrow \mathbb{E} [X]\geq \mathbb{E} [Y]$.

\item for all $X,Y$ in $B_{b}(\Omega)$, $\mathbb{E} [X+Y]\leq
\mathbb{E} [X]+\mathbb{E} [Y]$.

\item for all $\lambda \geq0$, $X\in B_{b}(\Omega)$,
 $\mathbb{E} [\lambda X]=\lambda \mathbb{E} [X]$.

\item for all $c\in \mathbb{R}$, $X\in
B_{b}(\Omega)$ , $\mathbb{E} [X+c]=\mathbb{E} [X]+c$.

\end{enumerate}
\end{theorem}

Moreover, it is also easy to check

\begin{theorem}
We have
\begin{enumerate}
\item Let $\mathbb{E}[X_n]$ and $\mathbb{E}[\sum _{n=1}^{\infty}X_n]$ be finite. Then
$\mathbb{E}[\sum _{n=1}^{\infty}X_n] \leq \sum _{n=1}^{\infty}
\mathbb{E}[X_n].$

\item Let $X_n \uparrow X$ and $\mathbb{E}[X_n]$, $\mathbb{E}[X]$ be finite. Then
$\mathbb{E}[X_n] \uparrow \mathbb{E}[X]$.
\end{enumerate}
\end{theorem}

\begin{definition}
The functional $\mathbb{E} [\cdot]$ is said to be regular if for each $\{X_{n}%
\}_{n=1}^{\infty}$ in $C_{b}(\Omega)$ such that $X_{n}\downarrow0$ \
on $\Omega$, we have $\mathbb{E} [X_{n}]\downarrow0$.
\end{definition}

Similar to Lemma \ref{Lemma1} we have:

\begin{theorem}
\label{Thm2}$\mathbb{\mathbb{E}}[\cdot]$ is regular if and only if
$\mathscr{P}$ is relatively compact.
\end{theorem}

\begin{proof}
\textquotedblleft$\Longrightarrow$\textquotedblright \ part: For
each sequence of closed subsets $F_{n}\downarrow \emptyset$ such
that $F_{n}$, $n=1,2,\cdots$, are non-empty (otherwise the proof is
trivial), there exists $\{g_{n}\}_{n=1}^{\infty}\subset
C_{b}(\Omega)$ satisfying
\[
0\leq g_{n}\leq1,\  \ g_{n}=1\text{ on }F_{n}\text{ and
}g_{n}=0\text{ on }\{ \omega \in \Omega:d(\omega,F_{n})\geq
\frac{1}{n}\}.
\]
We set $f_{n}=\wedge_{i=1}^{n}g_{i}$, it is clear that $f_{n}\in
C_{b}(\Omega)$ and $\mathbf{1}_{F_{n}}\leq f_{n}\downarrow0$. \
$\mathbb{E} [\cdot]$ is regular implies $\mathbb{E}
[f_{n}]\downarrow0$ and thus $c(F_{n})\downarrow 0$. It follows from
Lemma \ref{Lemma1} that $\mathscr{P}$ is relatively compact.
\newline \textquotedblleft$\Longleftarrow$\textquotedblright \ part:
For each $\left \{  X_{n}\right \}  _{n=1}^{\infty}\subset
C_{b}(\Omega)$ such
that $X_{n}\downarrow0$, we have%
\[
\mathbb{E} [X_{n}]=\sup_{P\in \mathscr{P}}E_{P}[X_{n}]=\sup_{P\in \mathscr{P}}%
\int_{0}^{\infty}P(\{X_{n}\geq t\})dt\leq
\int_{0}^{\infty}c(\{X_{n}\geq t\})dt.
\]
For each fixed $t>0$, $\{X_{n}\geq t\}$ is a closed subset and
$\{X_{n}\geq t\} \downarrow \emptyset$ as $n\uparrow \infty$. By
Lemma \ref{Lemma1}, $c(\{X_{n}\geq t\})\downarrow0$ and thus
$\int_{0}^{\infty}c(\{X_{n}\geq t\})dt\downarrow0$. Consequently
$\mathbb{E} [X_{n}]\downarrow0$.\end{proof}

%For the case where $\mathscr{P}$ is relatively compact, we denote by
%$\mathscr{\bar{P}}$ the closure of $\mathscr{P}$ in the topology of
%weak convergence. Then $\mathscr{\bar{P}}$ is weakly compact. From
%the definition
%of weak convergence%
%\[
%\mathbb{E} ^{\mathscr{P}}[X]=\mathbb{E} ^{\mathscr{\bar{P}}}[X],\  \
%\ X\in C_{b}(\Omega).
%\]
%In the rest of this section we will always assume that $\mathscr{P}$
%is \textit{weakly compact}. But some results still hold for general
%case.
%Thanks to the capacitability theorem of
%Choquet (see \cite{Cho, deme, fedlp, HuSt}), we have the following
%result which states that $c$ is \textit{regular} (with the
%capacity-related vocabulary) in the sense that it satisfies the
%Prokhorov's property. The proof is given in Appendix.

\subsection{Functional spaces}

We set, for $p>0$,

\begin{itemize}
\item $\mathcal{L}^{p}:=\{X\in L^{0}(\Omega):\mathbb{E} [|X|^{p}]=\sup
_{P\in \mathscr{P}}E_{P}[|X|^{p}]<\infty \}$;\  \

\item $\mathcal{N}^{p}:=\{X\in L^{0}(\Omega):\mathbb{E}
[|X|^{p}]=0\}$;

\item $\mathcal{N}:=\{X\in L^{0}(\Omega):X=0$, $c$-q.s.$\}$.
\end{itemize}

It is seen that $\mathcal{L}^{p}$ and $\mathcal{N}^{p}$ are linear
spaces and $\mathcal{N}^{p}=\mathcal{N}$, for each $p>0$.

We denote $\mathbb{L}^{p}:=\mathcal{L}^{p}/\mathcal{N}$. As usual,
we do not take care about the distinction between classes and their
representatives.\newline

\begin{lemma}
{\label{markov}} Let $X\in \mathbb{L}^{p}$. Then for each $\alpha>0$
\[
c(\{| X | >\alpha \})\leq \displaystyle
\frac{\mathbb{E}[|X|^{p}]}{\alpha^{p}}.
\]

\end{lemma}

\begin{proof} Just apply Markov inequality under each $P\in
\mathscr{P}$.
\end{proof}

Similar to the classical results, we get the following proposition
and the proof is omitted which is similar to the classical
arguments.
\begin{proposition}
\label{Prop3}We have
\begin{enumerate}
\item For each $p\geq1$, $\mathbb{L}^{p}$ is a Banach
space under the norm $\left \Vert X\right \Vert _{p}:=\left(  \mathbb{{E}%
}[|X|^{p}]\right)  ^{\frac{1}{p}}$.

\item For each $p<1$,
$\mathbb{L}^{p}$
is a complete metric space under the distance \newline $d(X,Y):= \mathbb{{E}%
}[|X-Y|^{p}]  $.
\end{enumerate}
\end{proposition}

\  \

We set%
\begin{align*}
\mathcal{L}^{\infty}  &  :=\{X\in L^{0}(\Omega):\exists \text{ a
constant
}M\text{, s.t. }|X|\leq M,\  \  \text{q.s.}\}; \\
\mathbb{L}^{\infty}  &  :=\mathcal{L}^{\infty}/\mathcal{N}.
\end{align*}

\begin{proposition}
\label{Prop4}Under the norm
\[
\left \Vert X\right \Vert _{\infty}:=\inf \left \{  M\geq0:|X|\leq
M,\  \  \text{q.s.}\right \}  ,
\]
$\mathbb{L}^{\infty}$ is a Banach space.
\end{proposition}

\begin{proof}From $\left \{  \left \vert X\right \vert >\left \Vert X\right
\Vert _{\infty}\right \}  =\cup_{n=1}^{\infty}\left \{ \left \vert
X\right \vert \geq \left \Vert X\right \Vert
_{\infty}+\frac{1}{n}\right \}  $ we know that $\left \vert X\right
\vert \leq \left \Vert X\right \Vert _{\infty}$, q.s., then it is
easy to check that $\left \Vert \cdot \right \Vert_{\infty} $ is a
norm. The proof of the completeness of $\mathbb{L}^{\infty}$ is
similar to the classical result. \end{proof}

With respect to the distance defined on $\mathbb{L}^{p}$, $p>0$, we
denote by

\begin{itemize}
\item $\mathbb{L}^{p}_{b}$ the completion of $B_{b}(\Omega)$.

\item $\mathbb{L}^{p}_{c}$ the completion of $C_{b}(\Omega)$.
\end{itemize}

By Proposition \ref{Prop3}, we have%
\[
\mathbb{L}^{p}_{c}\subset \mathbb{L}^{p}_{b}\subset \mathbb{L}^{p},\
\  \ p>0.
\]

The following Proposition is obvious and the proof is left to the
reader.

\begin{proposition}
We have
\begin{enumerate}
\item Let $p,q>1$, $\frac{1}{p}+\frac{1}{q}=1$. Then $X\in \mathbb{L}^{p}$ and
$Y\in \mathbb{L}^{q}$ implies
\[
XY\in \mathbb{L}^{1}\text{ and }\mathbb{E} [|XY|]\leq \left(
\mathbb{E} [|X|^{p}]\right)  ^{\frac{1}{p}}\left(
\mathbb{E}[|Y|^{q}]\right)
^{\frac{1}{q}};%
\]
\newline Moreover $X\in \mathbb{L}^{p}_{c}$ and $Y\in \mathbb{L}^{q}_{c}$  implies
$XY\in \mathbb{L}^{1}_{c}$;

\item $\mathbb{L}^{p_{1}}\subset \mathbb{L}^{p_{2}}$,\ $\mathbb{L}^{p_{1}}%
_{b}\subset \mathbb{L}^{p_{2}}_{b}$, $\mathbb{L}^{p_{1}}_{c}\subset
\mathbb{L}^{p_{2}}_{c}$, $0<p_{2}\leq p_{1}\leq \infty$;

\item $\left \Vert X\right \Vert _{p}\uparrow \left \Vert X\right \Vert _{\infty}$,
for each $X\in \mathbb{L}^{\infty}$.

%\item For $p\geq1$ and $X,Y\in \mathbb{L}^{p}$ we have
%\[
%\left(  \mathbb{E} [|X+Y|^{p}]\right)  ^{\frac{1}{p}}\leq \left(
%\mathbb{E} [|X|^{p}]\right)  ^{\frac{1}{p}}+\left(
%\mathbb{E}[|Y|^{p}]\right) ^{\frac{1}{p}}.
%\]

\end{enumerate}
\end{proposition}

\begin{proposition}
Let $p\in(0,\infty]$ and $(X_{n} )$ be a sequence in
$\mathbb{L}^{p}$ which converges to $X$ in $\mathbb{L}^{p}$. Then
there exists a subsequence $(X_{n_{k}})$ which converges to $X$
quasi-surely in the sense that it converges to $X$ outside a polar
set.
\end{proposition}

\begin{proof}Let us assume $p\in (0,\infty)$, the case
$p=\infty$ is obvious since the convergence in $\linf$ implies the
convergence in $\lp$ for all $p$.\\
One can extract a subsequence $(X_{n_k})$ such that
\[ \E [|X-X_{n_k}|^p]\leq 1/k^{p+2},\  \  \ k\in \N .\]
We set for all $k$
\[ A_k =\{ |X-X_{n_k}|>1/k \},\]
then as a consequence of the Markov property (Lemma \ref{markov})
and the Borel-Cantelli Lemma \ref{BorelC},
$c(\overline{\lim}_{k\rightarrow \infty}A_k )=0$. As it is clear
that on $(\overline{\lim}_{k\rightarrow \infty}A_k )^c$, $(X_{n_k})$
converges to $X$, the proposition is proved. \end{proof}

We now give a description of $\mathbb{L}^{p}_{b}$.

\begin{proposition}
\label{Prop5}For each $p>0$,%
\[
\mathbb{L}^{p}_{b}=\{X\in \mathbb{L}^{p}:\lim_{n\rightarrow
\infty}\mathbb{E} [|X|^{p}\mathbf{1}_{\{|X|>n\}}]=0\}.
\]

\end{proposition}

\begin{proof}We denote $J_{p}=\{X\in
\mathbb{L}^{p}:\lim_{n\rightarrow \infty}\E
[|X|^{p}\mathbf{1}_{\{|X|>n\}}]=0\}$. For each $X\in
J_{p}$ let $X_{n}=(X\wedge n)\vee(-n)\in B_{b}(\Omega)$. We have%
\[
\E [|X-X_{n}|^{p}]\leq \E [|X|^{p}\mathbf{1}%
_{\{|X|>n\}}]\rightarrow0\text{, as }n\rightarrow \infty \text{. }%
\]
Thus $X\in \mathbb{L}^{p}_{b}$.

On the other hand, for each $X\in \mathbb{L}^{p}_{b}$, we can find a
sequence $\left \{ Y_{n}\right \}
_{n=1}^{\infty}$ in $B_{b}(\Omega)$ such that $\E [|X-Y_{n}%
|^p]\rightarrow0$. Let $y_{n}=\sup_{\omega \in
\Omega}|Y_{n}(\omega)|$ and $X_{n}=(X\wedge y_{n})\vee (-y_{n})$.
Since $|X-X_{n}|\leq|X-Y_{n}|$, we have $\E [|X-X_{n}|^p
]\rightarrow0$. This clearly implies that for any sequence
$(\alpha_n )$
tending to $\infty$, $\lim_{n\rightarrow \infty} \E [|X-(X\wedge \alpha_{n})\vee (-\alpha_{n})|^p ]=0$.\\
Now we have, for all $n\in \N$,
\begin{align*}
\E [|X|^{p}\mathbf{1}_{\{|X|>n\}}]  &  =\mathbb{{E}%
}[(|X|-n+n)^{p}\mathbf{1}_{\{|X|>n\}}]\\
&  \leq(1\vee2^{p-1})\left( \E [(|X|-n)^{p}\mathbf{1}_{\{|X|>n\}}%
]+n^{p}c(|X|>n)\right).
\end{align*}
The first term of the right hand side tends to $0$ since%
\[
\E [(|X|-n)^{p}\mathbf{1}_{\{|X|>n\}}]= \E %
[|X-(X\wedge {n})\vee {(-n)}|^{p}]\rightarrow0.
\]
For the second term, since%
\[
\frac{n^{p}}{2^{p}}\mathbf{1}_{\{|X|>n\}}\leq(|X|-\frac{n}{2})^{p}%
\mathbf{1}_{\{|X|>n\}}\leq(|X|-\frac{n}{2})^{p}\mathbf{1}_{\{|X|>\frac{n}%
{2}\}},
\]
we have
\[
\frac{n^{p}}{2^{p}}c(|X|>n)=\frac{n^{p}}{2^{p}}\E [\mathbf{1}%
_{\{|X|>n\}}]\leq \mathbb{E}[(|X|-\frac{n}{2}%
)^{p}\mathbf{1}_{\{|X|>\frac{n}{2}\}}]\rightarrow0.
\]
Consequently $X\in J_{p}$.
\end{proof}

\begin{proposition}
\label{Prop12}Let $X\in \mathbb{L}^{1}_{b}$. Then for each
$\varepsilon>0$, there exists a $\delta>0$, such that for all $A\in
\mathscr{B}(\Omega)$ with $c(A)\leq \delta$, we have
$\mathbb{E}[|X|\mathbf{1}_{A}]\leq \varepsilon$.
\end{proposition}

\begin{proof}For each $\varepsilon>0$, by Proposition \ref{Prop5}, there exists an $N>0$ such that
$\E [|X|\mathbf{1}_{\{|X|>N\}}]\leq \frac{\varepsilon}{2}$. Take $\delta=\frac{\varepsilon}{2N}$. Then for a subset $A\in \mathscr{B}%
(\Omega)$ with $c(A)\leq \delta$, we have%
\begin{align*}
\E [|X|\mathbf{1}_{A}]  &  \leq \E [|X|\mathbf{1}%
_{A}\mathbf{1}_{\{|X|>N\}}]+\E [|X|\mathbf{1}_{A}\mathbf{1}%
_{\{|X|\leq N\}}]\\
&  \leq \E [|X|\mathbf{1}_{\{|X|>N\}}]+Nc(A)\leq \varepsilon
\text{.}%
\end{align*}
\end{proof}

It is important to note that not every element in $\mathbb{L}^{p}$
satisfies the condition $\lim_{n\rightarrow \infty}\mathbb{E}
[|X|^{p}\mathbf{1}_{\{|X|>n\}}]=0$. We give the following two
counterexamples to show that $\mathbb{L}^{1}$ and
$\mathbb{L}^{1}_{b}$ are different spaces even under the case that
$\mathscr{P}$ is weakly compact.

\begin{example}
\label{Exm2}Let $\Omega=\mathbb{N}$, $\mathscr{P}=\{P_{n}:n\in
\mathbb{N}\}$ where $P_{1}(\{1\})=1$ and
$P_{n}(\{1\})=1-\frac{1}{n}$, $P_{n}(\{n\})=\frac{1}{n}$, for
$n=2,3,\cdots$. $\mathscr{P}$ is weakly compact. We consider a
function $X$ on $\mathbb{N}$ defined by $X(n)=n$, $n\in \mathbb{N}$.
We have $\mathbb{E} [|X|]=2$ but $\mathbb{E}
[|X|\mathbf{1}_{\{|X|>n\}}]=1\not \rightarrow 0$. In this case,
$X\in \mathbb{L}^{1}$ but $X\not \in \mathbb{L}^{1}_{b}$.
\end{example}

\begin{example}
\label{Exm3}Let $\Omega=\mathbb{N}$, $\mathscr{P}=\{P_{n}:n\in
\mathbb{N}\}$ where $P_{1}(\{1\})=1$ and
$P_{n}(\{1\})=1-\frac{1}{n^{2}}$,  $P_{n}(\{kn\})=\frac{1}{n^{3}}$,
$k=1,2,\ldots,n$,for $n=2,3,\cdots$. $\mathscr{P}$ is weakly
compact. We consider a function $X$ on $\mathbb{N}$
defined by $X(n)=n$, $n\in \mathbb{N}$. We have $\mathbb{E}%
[|X|]=\frac{25}{16}$ and $n\mathbb{E}[\mathbf{1}_{\{|X|\geq
n\}}]=\frac {1}{n}\rightarrow0$, but
$\mathbb{E}[|X|\mathbf{1}_{\{|X|\geq n\}}]=\frac
{1}{2}+\frac{1}{2n}\not \rightarrow 0$. In this case, $X$ is in $\mathbb{L}%
^{1}$, continuous and $n\mathbb{E}[\mathbf{1}_{\{|X|\geq
n\}}]\rightarrow0$, but it is not in $\mathbb{L}^{1}_{b}$.
\end{example}

\subsection{Properties of elements in $\mathbb{L}^{p}_{c}$}

\begin{definition}
A mapping $X$ on $\Omega$ with values in a topological space is said
to be \textrm{quasi-continuous} (q.c.) if
\[
\forall \varepsilon>0 ,\; \makebox{there exists an open set } O\
\makebox{ with}\ c(O)<\varepsilon \makebox{ such that } X|_{O^{c}}
\makebox{ is continuous}.
\]

\end{definition}

\begin{definition}
We say that $X:\Omega \rightarrow \mathbb{R}$ has a quasi-continuous
version if there exists a quasi-continuous function $Y:\Omega
\rightarrow \mathbb{R}$ with $X=Y$ q.s..
\end{definition}

%\begin{remark}
%%In the rest of this paper, we say that $X\in \mathbb{L}^{p}$ has a
%%quasi-continuous version which means that there exists a
%%quasi-continuous function $Y\in L^0(\Omega)$ such that $X=Y$ q.s..
%Let $X$ be quasi-continuous and $Y=X$ q.s.. If $Y$ is
%quasi-continuous which is still an open problem.
%\end{remark}
%\begin{lemma}
%\label{Lemma7}If $X$ is quasi-continuous and $Y=X$ q.s., then $Y$ is
%also quasi-continuous.
%\end{lemma}
%
%\begin{proof}For each $\varepsilon>0$, since $c(\{X\not =Y\})=0$,
% there exists an open subset $G_{1}\supseteq \{X\not =Y\}$ such
%that $c(G_{1})<\frac{\varepsilon}{2}$. We can also find an open set
%$G_{2}$ with $c(G_{2} )<\frac{\varepsilon}{2}$ such that $X$ is
%continuous on $G_{2}^{c}$. It follows that the open subset
%$G=G_{1}\cup G_{2}$ satisfies $c(G)<\varepsilon$ and $Y$ is
%continuous on $G^{c}$.\end{proof}
\begin{proposition}
{\label{qc}} Let $p>0$. Then each element in $\mathbb{L}^{p}_{c}$
has a quasi-continuous version.
\end{proposition}

\begin{proof} Let $(X_n )$ be a Cauchy sequence in $C_b (\Omega)$ for the distance on $\lp$.
Let us choose a subsequence $(X_{n_k} )_{k\geq 1}$ such that
\[ \E [|X_{n_{k+1}}-X_{n_k}|^p]\leq 2^{-2k},\  \  \  \forall k\geq 1,\]
and set for all $k$,
\[ A_k =\bigcup_{i=k}^{\infty} \{ |X_{n_{i+1}}-X_{n_i}|> 2^{-i/p}\}.\]
Thanks to the subadditivity property and the Markov inequality, we
have
\[ c(A_k )\leq \sum_{i=k}^{\infty} c(|X_{n_{i+1}}-X_{n_i}|>
2^{-i/p})\leq \sum_{i=k}^{\infty} 2^{-i}=2^{-k+1}.\] As a
consequence, $\lim_{k\rightarrow \infty} c(A_k )=0$, so the Borel
set $A=\bigcap_{k=1}^{\infty}A_{k}$ is polar.\\
As each $X_{n_k}$ is continuous, for all $k\geq 1$, $A_k$ is an open
set. Moreover, for all $k$, $(X_{n_i})$ converges uniformly on
$A_k^c$ so that the limit is continuous on each $A_k ^c$. This
yields the result.
\end{proof}

The following theorem gives a concrete characterization of the space
$\mathbb{L}^{p}_{c}$.

\begin{theorem}
\label{Thm8}For each $p>0$,%
\[
\mathbb{L}^{p}_{c}=\{X\in \mathbb{L}^{p} :X\text{ has a quasi-continuous version, }%
\lim_{n\rightarrow \infty}\mathbb{\mathbb{E} }[|X|^{p}\mathbf{1}_{\{|X|>n\}}%
]=0\}.
\]

\end{theorem}

\begin{proof}We denote
\[
J_{p}=\{X\in \mathbb{L}^{p}:X\text{ has a quasi-continuous version,
}\lim_{n\rightarrow \infty}\E [|X|^{p}\mathbf{1}_{\{|X|>n\}}]=0\}.
\]
Let $X\in \lpc$, we know by Proposition \ref{qc} that $X$ has a
quasi-continuous version. Since $X\in \mathbb{L}^{p}_b$,
we have by Proposition \ref{Prop5} that $\lim_{n\rightarrow \infty}%
\E [|X|^{p}\mathbf{1}_{\{|X|>n\}}]=0$. Thus $X\in J_{p}$.
\newline On the other hand, let $X\in J_{p}$ be quasi-continuous. Define $Y_{n}=(X\wedge n)\vee(-n)$ for all
$n\in \N$. As $\E [|X|^{p}\mathbf{1}_{\{|X|>n\}}]\rightarrow0$, we
have $\E [|X-Y_{n}|^{p}]\rightarrow0$. \\
Moreover, for all $n\in \N$, as  $Y_{n}$ is quasi-continuous , there
exists a closed set $F_{n}$ such that
$c(F_{n}^{c})<\frac{1}{n^{p+1}}$ and $Y_{n}$ is continuous on
$F_{n}$. It follows from Tietze's extension theorem that there
exists $Z_{n}\in C_{b}(\Omega)$ such that
\[
|Z_{n}|\leq n\text{ and }Z_{n}=Y_{n}\text{ \ on }F_{n}.
\]
We then have%
\[
\E [|Y_{n}-Z_{n}|^{p}]\leq(2n)^{p}c(F_{n}^{c})\leq \frac{(2n)^{p}%
}{n^{p+1}}.
\]
So $\mathbb{E}[|X-Z_{n}|^{p}]\leq(1\vee2^{p-1})(\mathbb{E}%
[|X-Y_{n}|^{p}]+\mathbb{E}[|Y_{n}-Z_{n}|^{p}])\ \rightarrow0,$ and
$X\in \lpc$.
\end{proof}

We give the following example to show that $\mathbb{L}^{p}_{c}$ is
different from $\mathbb{L}^{p}_{b}$ even under the case that
$\mathscr{P}$ is weakly compact.

\begin{example}
\label{Exm1}Let $\Omega=[0,1]$, $\mathscr{P}=\{ \delta_{x}:x\in
\lbrack0,1]\}$ is weakly compact. It is seen that
$\mathbb{L}^{p}_{c}=C_{b}(\Omega)$ which is different from
$\mathbb{L}^{p}_{b}$.
\end{example}

We denote $\mathbb{L}^{\infty}_{c}:=\{X\in \mathbb{L}^{\infty}:X$
has a quasi-continuous version$\}$, we have

\begin{proposition}
\label{Prop9}$\mathbb{L}^{\infty}_{c}$ is a closed linear subspace
of $\mathbb{L}^{\infty}$.
\end{proposition}

\begin{proof} For each Cauchy sequence $\left \{  X_{n}\right \}
_{n=1}^{\infty}$ of $\mathbb{L}^{\infty}_{c}$ under $\left \Vert
\cdot \right \Vert _{\infty}$, we can find a subsequence $\left \{  X_{n_{i}%
}\right \}  _{i=1}^{\infty}$ such that $\left \Vert X_{n_{i+1}}-X_{n_{i}%
}\right \Vert _{\infty}\leq 2^{-i}$. We may further assume that each
$X_n$ is quasi-continuous. Then it is easy to prove that for each
$\varepsilon>0$, there exists an open set $G$ such that
$c(G)<\varepsilon$ and $\left \vert X_{n_{i+1}}-X_{n_{i}}\right
\vert \leq 2^{-i}$ for all $i\geq 1$ on $G^{c}$, which implies that
the limit belongs to $\mathbb{L}^{\infty}_{c}$.
%Thus there exists an
%$A\in \mathscr{B}(\Omega)$ with $c(A)=0$ such that $\left \vert X_{n_{i+1}%
%}-X_{n_{i}}\right \vert \leq 2^{-i}$ on $A^{c}$ for all $i\geq 1$.
%We define a function $X$ by letting $X=0$ on $A$, and
%$X=\lim_{i\rightarrow \infty}X_{n_{i}}$ on $A^{c}$. It is easy to
%prove that  $X\in \mathbb{L}^{\infty}_{c}$ and $\left \Vert X_{n}-X
%\right \Vert _{\infty} \rightarrow 0$.
\end{proof}

As an application of Theorem \ref{Thm8}, we can easily get the
following results.

\begin{proposition}
\label{Prop10}Assume that $X:\Omega \rightarrow \mathbb{R}$ has a
quasi-continuous version and that there exists a function $f:\mathbb{R}^{+}%
\rightarrow \mathbb{R}^{+}$ satisfying $\lim_{t\rightarrow
\infty}\frac {f(t)}{t^{p}}=\infty$ and $\mathbb{E} [f(|X|)]<\infty$.
Then $X\in \mathbb{L}^{p}_{c}$.
\end{proposition}

\begin{proof} For each $\varepsilon>0$, there exists an $N>0$ such that
$\frac{f(t)}{t^{p}}\geq \frac{1}{\varepsilon}$, for all $t\geq N$. Thus%
\[
\E [|X|^{p}\mathbf{1}_{\{|X|>N\}}]\leq \varepsilon \mathbb{E}%
[f(|X|)\mathbf{1}_{\{|X|>N\}}]\leq \varepsilon \E [f(|X|)]\text{.}%
\]
Hence $\lim_{N\rightarrow \infty}\E [|X|^{p}\mathbf{1}_{\{|X|>N\}}%
]=0$. From Theorem \ref{Thm8} we infer $X\in$ $\lpc$. \end{proof}

\begin{lemma}
\label{Lemma14}Let $\left \{  P_{n}\right \}  _{n=1}^{\infty}\subset
\mathscr{P}$ converge weakly to $P\in \mathscr{P}$. Then for each
$X\in \mathbb{L}^{1}_{c}$, we have $E_{P_{n}}[X]\rightarrow
E_{P}[X]$.
\end{lemma}

\begin{proof} We may assume that $X$ is quasi-continuous, otherwise we can
consider its quasi-continuous version which does not change the
value $E_{Q}$ for each $Q \in \mathscr{P}$. For each
$\varepsilon>0$, there exists an $N>0$ such that $\E
[|X|\mathbf{1}_{\left \{ |X|>N\right \} }]<\frac{\varepsilon }{2}$.
Set $X_{N}=(X\wedge N)\vee (-N)$. We can find an open subset $G$
such that $c(G)<\frac{\varepsilon}{4N}$ and $X_{N}$ is continuous on
$G^{c}$. By Tietze's extension theorem, there exists $Y\in
C_{b}(\Omega)$ such that $\left \vert Y\right \vert \leq N$ and
$Y=X_{N}$ on $G^{c}$. Obviously, for each $Q \in \mathscr{P}$,
\begin{align*}
|E_{Q}[X]-E_{Q}[Y]|  & \leq E_{Q}[|X-X_{N}|]+E_{Q}[|X_{N}-Y|]\\
& \leq \frac{\varepsilon}{2}+2N\frac{\varepsilon}{4N}=\varepsilon \text{. }%
\end{align*}
It then follows that
\[
\limsup_{n\rightarrow \infty}E_{P_{n}}[X]\leq \lim_{n\rightarrow \infty}E_{P_{n}%
}[Y]+\varepsilon=E_{P}[Y]+\varepsilon \leq E_{P}[X]+2\varepsilon,
\]
and similarly $\liminf_{n\rightarrow \infty}E_{P_{n}}[X]\geq E_{P}%
[X]-2\varepsilon$. Since $\varepsilon$ can be arbitrarily small, we
then have $E_{P_{n}}[X]\rightarrow E_{P}[X]$.
\end{proof}

\begin{remark}
For continuous $X$, the above lemma is Lemma 3.8.7 in \cite{bog}.
\end{remark}

Now we give an extension of Theorem \ref{Thm2}.
\begin{theorem}
\label{Thm15}Let $\mathscr{P}$ be weakly compact and let $\left \{
X_{n}\right \}  _{n=1}^{\infty}\subset \mathbb{L}
^{1}_{c}$ be such that $X_{n}\downarrow X$, q.s.. Then $\mathbb{E}%
[X_{n}]\downarrow \mathbb{E}[X]$.
\end{theorem}

\begin{remark}
It is important to note that $X$ does not necessarily belong to $\mathbb{L}%
^{1}_{c}$.
\end{remark}

\begin{proof}For the case $\E [X]>-\infty$, if there exists a
$\delta>0$ such that $\E [X_{n}]>\E [X]+\delta$,
$n=1,2,\cdots$, we then can find a $P_{n}\in \mathscr{P}$ such that $E_{P_{n}%
}[X_{n}]>\E [X]+\delta-\frac{1}{n}$, $n=1,2,\cdots$. Since
$\mathscr{P}$ is weakly compact, we then can find a subsequence
$\left \{ P_{n_{i}}\right \}  _{i=1}^{\infty}$ that converges weakly
to some $P\in \mathscr{P}$. From which it follows that
\begin{align*}
E_{P}[X_{n_{i}}]  & =\lim_{j\rightarrow \infty}E_{P_{n_{j}}}[X_{n_{i}}%
]\geq \limsup_{j\rightarrow \infty}E_{P_{n_{j}}}[X_{n_{j}}]
\\&\geq \limsup _{j\rightarrow \infty}\{ \E
[X]+\delta-\frac{1}{n_{j}}\} =\E [X]+\delta \text{,\  \
}i=1,2,\cdots.
\end{align*}
Thus   $E_{P}[X]\geq \E [X]+\delta$. This contradicts  the
definition of $\E [\cdot]$. The proof for the case $\E [X]=-\infty$
is analogous.\end{proof}

We immediately have the following corollary.
\begin{corollary}\label{newco33}
Let $\mathscr{P}$ be weakly compact and let $\left \{  X_{n}\right
\} _{n=1}^{\infty}$ be a sequence in $\mathbb{L}^{1}_{c}$
decreasingly converging to $0$ q.s.. Then $\mathbb{E}
[X_{n}]\downarrow0$.
\end{corollary}

\subsection{Kolmogorov's criterion}

\begin{definition}
Let $I$ be a set of indices, $(X_{t})_{t\in I}$ and $(Y_{t} )_{t\in
I}$ be two processes indexed by $I$ . We say that $Y$ is \textrm{a
quasi-modification} of $X$ if for all $t\in I$, $X_{t} =Y_{t} $
q.s..
\end{definition}

\begin{remark}
In the above definition, quasi-modification is also called
modification in some papers.
\end{remark}

We now give a Kolmogorov criterion for a process indexed by
$\mathbb{R}^{d}$ with $d\in \mathbb{N}$.

\begin{theorem}
Let $p>0$ and $(X_{t} )_{t\in[0,1]^{d}}$ be a process such that for
all $t\in[0,1]^{d}$, $X_{t}$ belongs to $\mathbb{L}^{p}$ . Assume
that there exist positive constants $c$ and $\varepsilon$ such that
\[
\mathbb{E} [|X_{t} -X_{s} |^{p}]\leq c|t-s|^{d+\varepsilon}.
\]
Then $X$ admits a modification $\tilde{X}$ such that
\[
\mathbb{E}\left[ \left(  \sup_{s\neq t} \displaystyle
\frac{|\tilde{X}_{t} -\tilde{X}_{s}|}{|t-s|^{\alpha}}\right)
^{p}\right] <\infty,
\]
for every $\alpha \in[0,\varepsilon/p)$. As a consequence, paths of
$\tilde{X}$ are quasi-surely Hölder continuous of order $\alpha$ for
every $\alpha<\varepsilon/p$ in the sense that there exists a Borel
set $N$ of capacity $0$ such that for all $w\in N^{c}$, the map
$t\rightarrow \tilde{X} (w)$ is Hölder continuous of order $\alpha$
for every $\alpha <\varepsilon/p$. Moreover, if
$X_t\in\mathbb{L}_c^p$ for each $t$, then we also have
$\tilde{X}_t\in\mathbb{L}_c^p$.
\end{theorem}

\begin{proof} Let $D$ be the set of dyadic points in $[0,1]^d$:
\[ D=\left \{ (\frac{i_1}{2^n} ,\cdots ,\frac{i_d}{2^n}); \ n\in \N , i_1
,\cdots ,i_d \in \{0,1,\cdots ,2^n \} \right \}.\] Let $\alpha \in
[0,\vep/p)$. We set
\[ M=\sup_{s,t\in D, s\neq t} \displaystyle \frac{|X_t
-X_s|}{|t-s|^{\alpha}}.\] Thanks to the classical Kolmogorov's
criterion (see Revuz-Yor \cite{ReYo}), we know that for any $P\in
\mathscr{P}$, $E_P [M^{p}]$ is finite and uniformly bounded with
respect to $P$ so that
\[ \E [M^{p}]=\sup_{P\in \PP}E_P
[M^{p}]<\infty.\] As a consequence, the map $t\rightarrow X_t $ is
uniformly continuous on $D$ quasi-surely and so we can define
\[ \forall t\in [0,1]^d ,\  \tX_t =\lim_{s\rightarrow t , s\in D} X_s
.\] It is now clear that $\tX$ satisfies the enounced properties.
\end{proof}

\section{$G$-Brownian motion under $G$-expectations}

In this section we consider the following path spaces:
$\Omega=C_{0}^{d}(\mathbb{R}^{+})$ the space of all
$\mathbb{R}^{d}$-valued continuous paths
$(\omega_{t})_{t\in\mathbb{R}^{+}}$, with $\omega_{0}=0$, equipped
with the distance
\[
\rho(\omega^{1},\omega^{2}):=\sum_{i=1}^{\infty}2^{-i}[(\max_{t\in
[0,i ]}|\omega_{t}^{1}-\omega_{t}^{2}|)\wedge1].
\]
It is clear that $(\Omega,\rho)$ is a complete separable metric
space. We also denote $\Omega_T=\{ \omega_{.\wedge T}:\omega \in
\Omega \}$ for each fixed $T \in [0,\infty)$.

Let $\mathcal{H}$ be a vector lattice of real functions defined on
$\Omega$ such that if $X_{1},\cdots,X_{n}\in \mathcal{H}$ then
$\varphi(X_{1},\cdots,X_{n})\in \mathcal{H}$ for each $\varphi \in
C_{b.Lip}(\mathbb{R}^{n})$, where $C_{b.Lip}(\mathbb{R}^{n})$
denotes the space of all bounded and Lipschitz functions on
$\mathbb{R}^{n}$.
\begin{definition}
A functional  $\mathbb{E}:\mathcal{H}\mapsto \mathbb{R}$ is called a
\textrm{sublinear expectation} on $\mathcal{H}$ if it satisfies:
\begin{enumerate}
\item \textbf{Monotonicity}: for all $X,Y$ in $\mathcal{H}$, $X\geq
Y\Longrightarrow \mathbb{E} [X]\geq \mathbb{E} [Y]$.

%%\item \textbf{Preservation of constants}: for all $c\in \mathbb{R}$,
%%$\mathbb{E} [c]=c$.

\item \textbf{Sub-additivity}: for all $X,Y$ in $\mathcal{H}$, $\mathbb{E}
[X+Y]\leq \mathbb{E} [X]+\mathbb{E} [Y]$.

\item \textbf{Positive homogeneity}: for all $\lambda \geq0$, $X\in \mathcal{H}%
$, $\mathbb{E} [\lambda X]=\lambda \mathbb{E} [X]$.

\item \textbf{Constant translatability}: for all $c\in \mathbb{R}$,
$X\in \mathcal{H}$, $\mathbb{E} [X+c]=\mathbb{E} [X]+c$.
\end{enumerate}
\label{def3.1}\end{definition}
%Let $\mathbb{E}[\cdot]:\mathcal{H}\mapsto \mathbb{R}$ be a sublinear
%expectation on $\mathcal{H}$ in the sense of Definition
%\ref{def3.1}.
A $d$-dimensional random vector $X$ with each
component in $\mathcal{H}$ is said to be \textit{$G$-normally
distributed} under the sublinear expectation $\mathbb{E}[\cdot]$ if
for each $\varphi \in
C_{b.Lip}(\mathbb{R}^{d})$, the function $u$ defined by%
\[
u(t,x):=\mathbb{E}[\varphi(x+\sqrt{t}X)],\ t\geq0,\  \ x\in
\mathbb{R}^{d}
\]
satisfies\ the following $G$-heat equation:%
\begin{align*}
\frac{\partial u}{\partial t}-G(D^{2}u)  &  =0,\  \  \text{on
}(t,x)\in
\lbrack0,\infty)\times \mathbb{R}^{d},\\
u(0,x)  &  =\varphi(x),
\end{align*}
where $D^{2}u$ is the Hessian matrix of $u$ , i.e.,
$D^{2}u=(\partial_{x_{i}x_{j}}^{2}u)_{i,j=1}^{d}$ and
\begin{equation}
G(A)=\frac{1}{2}\sup_{\gamma\in\Theta}\text{tr}[\gamma\gamma^{T}A],
\  \ A=(A_{ij})_{i,j=1}^{d}\in\mathbb{S}_{d}. \label{G}
\end{equation}
$\mathbb{S}_{d}$ denotes the space of $d\times d$ symmetric
matrices. $\Theta$ is a given non empty, bounded and closed subset
of $\mathbb{R}^{d\times d}$ which is the space of all $d\times d$
matrices.

\begin{remark} The above $G$-heat equation has a unique viscosity
solution. We refer to \cite{CIL} for the definition, existence,
uniqueness and comparison theory of this type of parabolic PDE (see
also \cite{Peng2007b} for our specific situation). If $G$ is
non-degenerate, i.e., there exists a $\beta>0$ such that
$G(A)-G(B)\geq \beta \textrm{Tr}[A-B]$ for each $A,\,B\in
\mathbb{S}_{d}$ with $A\geq B$, then the above $G$-heat equation has
a unique $C^{1,2}$-solution (see e.g. \cite{Wang}).
\end{remark}

We consider the canonical process: $B_t(\omega)=\omega_t$, $t\in
[0,\infty)$, for $\omega\in \Omega$. We introduce the space of
finite dimensional cylinder random variables: for each fixed
$T\geq0$, we set
\[
L_{ip}(\Omega_T):=\{ \varphi(B_{t_1},B_{t_2},\cdots,B_{t_n}):
\forall n\geq1, t_1,\cdots,t_n \in [0,T], \forall \varphi \in
C_{b.Lip}(\mathbb{R}^{d\times n}) \}, \]
 It is clear that
$L_{ip}(\Omega_t)\subseteq L_{ip}(\Omega_T)\subset C_b(\Omega_T)$,
for $t\leq T$. We also denote
\[
L_{ip}(\Omega):=\bigcup_{n=1}^{\infty}L_{ip}(\Omega_n)\subset
C_b(\Omega).
\]
We can construct (see \cite{Peng2006a,Peng2006b}) a consistent
sublinear expectation called $G-$expectation $\mathbb{E}[\cdot]$ on
$L_{ip}(\Omega)$, such that $B_1$ is $G-$normally distributed under
$\mathbb{E}[\cdot]$ and for each $s$, $t\geq0$ and
$t_{1},\cdots,t_{N}\in \lbrack0,t]$ we have
\begin{equation}
\mathbb{E}[\varphi(B_{t_{1}},\cdots,B_{t_{N}},B_{t+s}-B_{t})]=\mathbb{E}%
[\psi(B_{t_{1}},\cdots,B_{t_{N}})], \label{GBM}\end{equation}
where $\psi(x_{1},\cdots,x_{N})=\mathbb{E}[\varphi(x_{1}%
,\cdots,x_{N},\sqrt{s}B_1)]$. Under $G-$expectation
$\mathbb{E}[\cdot]$, the canonical process $\{ B_t : t \geq 0 \}$ is
called $G-$Brownian motion.
\begin{remark}
Relation (\ref{GBM}) implies that the increments of $B$ are
independent and stationary distributed with respect to the sublinear
expectation $\mathbb{E}[\cdot]$.  The condition that $B_1$ is
$G$-normally distributed can be also automatically obtained provided
that $\mathbb{E}[|B_t|^4]\leq Ct^2$ (see \cite{Peng2007b}).
\end{remark}

The topological completion of $L_{ip}(\Omega_T)$ (resp.
$L_{ip}(\Omega)$) under the Banach norm $\mathbb{E}[|\cdot|]$ is
denoted by $L_{G}^{1}(\Omega_T)$ (resp. $L_{G}^{1}(\Omega)$).
$\mathbb{E}[\cdot]$ can be extended uniquely to a sublinear
expectation on $L_{G}^{1}(\Omega)$.

In the previous section the sublinear expectation
$\mathbb{E}[\cdot]$ is induced as an upper expectation associated to
a family $\mathscr{P}$ of probability measures. In this Section
$\mathbb{E}[\cdot]$ will always be the $G-$expectation. We will
prove that $C_b(\Omega)\subset L_{G}^{1}(\Omega)$ and that, in fact,
$\mathbb{E}[\cdot]$ is the upper expectation of a weakly compact
family $\mathscr{P}$ on $\Omega$, thus all results in Section 2 hold
true.

\subsection{$G$-Expectation as an upper-Expectation}

In this subsection we will construct a family $\mathscr{P}$ of
probability measures on $\Omega$, for which the upper expectation
coincides with the $G$-expectation $\mathbb{E}[\cdot]$ on
$L_{ip}(\Omega)$.

Let $(\Omega,\mathcal{F},P)$ be a probability space and
$(W_{t})_{t\geq 0}=(W_{t}^{i})_{i=1,t\geq0}^{d}$ a $d$-dimensional
Brownian motion in this space. The filtration generated by $W$ is
denoted by
\[
\mathcal{F}_{t}:=\sigma \{W_{u},0\leq u\leq t\} \vee \mathcal{N},\  \  \mathbb{F}%
=\{ \mathcal{F}_{t}\}_{t\geq0},%
\]
where $\mathcal{N}$ is the collection of $P$--null subsets. We also
denote, for a fixed $s\geq0$,
\[
\mathcal{F}_{t}^{s}:=\sigma \{W_{s+u}-W_{s},0\leq u\leq t\} \vee \mathcal{N}%
,\  \  \mathbb{F}^{s}:=\{ \mathcal{F}_{t}^{s}\}_{t\geq0}.
\]
Let $\Theta$ be a given bounded and closed subset in
$\mathbb{R}^{d\times d}$. We denote by $\mathcal{A}_{t,T}^{\Theta}$,
the collection of all $\Theta $-valued $\mathbb{F}$-adapted process
on an interval $[t,T]\subset \lbrack0,\infty)$. For each fixed
$\theta \in \mathcal{A}_{t,T}^{\Theta}$ we denote
\[
B_{T}^{t,\theta}:=\int_{t}^{T}\theta_{s}dW_{s}.
\]
In this section we will prove that, for each $n=1,2,\cdots$,
$\varphi \in C_{b.Lip}(\mathbb{R}^{d\times n})$ and $0\leq
t_{1},\cdots,t_{n}<\infty$, the $G$-expectation defined in
\cite{Peng2006a, Peng2006b} can be equivalently defined by
\[
\mathbb{E}[\varphi(B_{t_{1}},B_{t_{2}}-B_{t_{1}},\cdots,B_{t_{n}}-B_{t_{n-1}%
})]=\sup_{\theta \in \mathcal{A}_{0,\infty}^{\Theta}}E_{P}[\varphi(B_{t_{1}%
}^{0,\theta},B_{t_{2}}^{t_{1},\theta},\cdots,B_{t_{n}}^{t_{n-1},\theta})].
\]

Given $\varphi \in C_{b.Lip}(\mathbb{R}^{n}\times \mathbb{R}^{d})$,
$0\leq t\leq T<\infty$ and $\zeta \in
L^{2}(\Omega,\mathcal{F}_{t},P;\mathbb{R}^{n})$, we define
\begin{equation}
\Lambda_{t,T}[\zeta]=\text{ess}\sup_{\theta \in \mathcal{A}_{t,T}^{\Theta}}%
E_{P}[\varphi(\zeta,\int_{t}^{T}\theta_{s}dW_{s})|\mathcal{F}_{t}].
\label{eq0}%
\end{equation}

\begin{lemma}
\label{lem1}For each $\theta^{1}$ and $\theta^{2}$ in $\mathcal{A}%
_{t,T}^{\Theta}$, there exists $\theta \in
\mathcal{A}_{t,T}^{\Theta}$ such that
\begin{equation}
E_{P}[\varphi(\zeta,B_{T}^{t,\theta})|\mathcal{F}_{t}]=E_{P}[\varphi
(\zeta,B_{T}^{t,\theta^{1}})|\mathcal{F}_{t}]\vee
E_{P}[\varphi(\zeta
,B_{T}^{t,\theta^{2}})|\mathcal{F}_{t}]. \label{eq1}%
\end{equation}
Consequently, there exists a sequence $\left \{  \theta^{i}\right \}
_{i=1}^{\infty}$ of $\mathcal{A}_{t,T}^{\Theta}$, such that
\begin{equation}
E_{P}[\varphi(\zeta,B_{T}^{t,\theta^{i}})|\mathcal{F}_{t}]\nearrow
\Lambda_{t,T}[\zeta],\  \ P\text{-a.s..} \label{eq2}%
\end{equation}
We also have, for each $s\leq t$,%
\begin{equation}
E_{P}[\text{ess}\sup_{\theta \in \mathcal{A}_{t,T}^{\Theta}}E_{P}%
[\varphi(\zeta,\int_{t}^{T}\theta_{r}dW_{r})|\mathcal{F}_{t}]|\mathcal{F}%
_{s}]=ess\sup_{\theta \in
\mathcal{A}_{t,T}^{\Theta}}E_{P}[\varphi(\zeta
,\int_{t}^{T}\theta_{r}dW_{r})|\mathcal{F}_{s}]. \label{eq3}%
\end{equation}

\end{lemma}

\begin{proof}
We set $A=\left \{  \omega:E_{P}[\varphi(\zeta,B_{T}^{t,\theta^{1}%
})|\mathcal{F}_{t}](\omega)\geq E_{P}[\varphi(\zeta,B_{T}^{t,\theta^{2}%
})|\mathcal{F}_{t}](\omega)\right \}  $ and take $\theta_{s}=I_{[t,T]}%
(s)(I_{A}\theta_{s}^{1}+I_{A^{C}}\theta_{s}^{2})$. Since%
\[
\varphi(\zeta,B_{T}^{t,\theta})=I_{A}\varphi(\zeta,B_{T}^{t,\theta^{1}%
})+I_{A^{C}}\varphi(\zeta,B_{T}^{t,\theta^{2}}),
\]
we derive (\ref{eq1})\ and then (\ref{eq2}). (\ref{eq3}) follows
from (\ref{eq1}) and Yan's commutation theorem (cf \cite{Yan} in
Chinese and Thm. a3 in the Appendix of \cite{Peng2004}).
\end{proof}

\begin{lemma}
\label{lem2}The mapping $\Lambda_{t,T}[\cdot]:L^{2}(\Omega,\mathcal{F}%
_{t},P;\mathbb{R}^{n})\rightarrow
L^{2}(\Omega,\mathcal{F}_{t},P;\mathbb{R})$ has the following
regularity properties: for each $\zeta$, $\zeta^{\prime}\in
L^{2}(\mathcal{F}_{t})$:\newline \textrm{(i)}
$\Lambda_{t,T}[\zeta]\leq
C_{\varphi}$.\newline \textrm{(ii)} $|\Lambda_{t,T}[\zeta]-\Lambda_{t,T}%
[\zeta^{\prime}]|\leq k_{\varphi}|\zeta-\zeta^{\prime}|.$\newline
where $C_{\varphi}=\sup_{(x,y)}\varphi(x,y)$ and $k_{\varphi}$ is
the Lipschitz constant of $\varphi$.
\end{lemma}

\begin{proof}
We only need to prove (ii). We have%
\begin{align*}
\Lambda_{t,T}[\zeta]-\Lambda_{t,T}[\zeta^{\prime}]  &  \leq \text{ess}%
\sup_{\mathcal{A}_{t,T}^{\Theta}}E_{P}[\varphi(\zeta,\int_{t}^{T}\theta
_{s}dW_{s})-\varphi(\zeta^{\prime},\int_{t}^{T}\theta_{s}dW_{s})|\mathcal{F}%
_{t}]\\
&  \leq k_{\varphi}|\zeta-\zeta^{\prime}|
\end{align*}
and, symmetrically,
$\Lambda_{t,T}[\zeta^{\prime}]-\Lambda_{t,T}[\zeta]\leq
k_{\varphi}|\zeta-\zeta^{\prime}|$. Thus (ii) follows.
\end{proof}

\begin{lemma}
For each $x\in \mathbb{R}^{n}$, $\Lambda_{t,T}[x]$ is a
deterministic function. Moreover,
\begin{equation}
\Lambda_{t,T}[x]=\Lambda_{0,T-t}[x]. \label{eq4}%
\end{equation}

\end{lemma}

\begin{proof}
Since the collection of processes $(\theta_{s})_{s\in \lbrack t,T]}$ with%
\[
\left \{  \theta_{s}=\sum_{j=1}^{N}I_{A_{j}}\theta_{s}^{j}:\left \{
A_{j}\right \}  _{j=1}^{N}\text{ is an
}\mathcal{F}_{t}\text{--partition of
}\Omega \text{, }\theta^{j}\in \mathcal{A}_{t,T}^{\Theta}\text{ is }%
(\mathbb{F}^{t})\text{--adapted}\right \}
\]
is dense in $\mathcal{A}_{t,T}^{\Theta}$, we can take a sequence
$\theta_{s}^{i}=\sum_{j=1}^{N_{i}}I_{A_{ij}}\theta_{s}^{ij}$ of this
type of
processes such that $E_{P}[\varphi(x,B_{T}^{t,\theta^{i}})|\mathcal{F}%
_{t}]\nearrow \Lambda_{t,T}[x]$. But%
\begin{align*}
E_{P}[\varphi(x,B_{T}^{t,\theta^{i}})|\mathcal{F}_{t}]  &  =\sum_{j=1}^{N_{i}%
}I_{A_{ij}}E_{P}[\varphi(x,B_{T}^{t,\theta^{ij}})|\mathcal{F}_{t}]=\sum_{j=1}^{N_{i}%
}I_{A_{ij}}E_{P}[\varphi(x,B_{T}^{t,\theta^{ij}})]\\
&  \leq \max_{1\leq j\leq N_{i}}E_{P}[\varphi(x,B_{T}^{t,\theta^{ij}})]=E_{P}%
[\varphi(x,B_{T}^{t,\theta^{ij_{i}}})],
\end{align*}
where, for each $i$, $j_{i}$ is a maximizer of $\left \{
E_{P}[\varphi (x,B_{T}^{t,\theta^{ij}})]\right \} _{j=1}^{N_{i}}$.
This implies that
\[
\lim_{i\rightarrow \infty}E_{P}[\varphi(x,B_{T}^{t,\theta^{ij_{i}}}%
)]=\Lambda_{t,T}[x],\  \  \  \text{a.s.}%
\]
and thus $\Lambda_{t,T}[x]$ is a deterministic number. In the above
proof, we know that
\[
ess\sup_{\theta \in
\mathcal{A}_{t,T}^{\Theta}}E_{P}[\varphi(x,B_{T}^{t,\theta
})|\mathcal{F}_{t}]=ess\sup_{\theta \in
\mathcal{A}_{0,T-t}^{t,\Theta}}E_{P}[\varphi(x,\int
_{0}^{T-t}\theta_{s}dW_{s}^{t})],
\]
where\ $W_{s}^{t}=W_{t+s}-W_{t}$, $s\geq0$, and $\mathcal{A}_{0,T-t}%
^{t,\Theta}$ is the collection of $\Theta$--valued and $\mathbb{F}^{t}%
$--adapted processes on $[0,T-t]$. Thus (\ref{eq4}) follows.
\end{proof}

We will denote $u_{t,T}(x):=\Lambda_{t,T}[x]$, $t\leq T$. By Lemma
\ref{lem2}, $u_{t,T}(\cdot)$ is a bounded and Lipschitz function.

\begin{lemma}
\label{lem3} For each $\zeta \in
L^{2}(\Omega,\mathcal{F}_{t},P;\mathbb{R}^{n})$, we have
\[
u_{t,T}(\zeta)=\Lambda_{t,T}[\zeta],\  \  \  \text{a.s.}.%
\]

\end{lemma}

\begin{proof}
By the above regularities of $\Lambda_{t,T}[\cdot]$ and
$u_{t,T}(\cdot)$ we only need to check the situation where $\zeta$
is a step function, i.e., $\zeta=\sum_{j=1}^{N}I_{A_{j}}x_{j}$,
where $x_{j}\in \mathbb{R}^{n}$ and $\left \{  A_{j}\right \}
_{j=1}^{N}$ is an $\mathcal{F}_{t}$--partition of of $\Omega$. For
each $x_{j}$, let $\left \{  \theta^{ij}\right \}  _{i=1}^{\infty }$
of $\mathcal{A}_{t,T}^{\Theta}$ be $(\mathcal{F}_{s}^{t})$--adapted
process such that
\[
\lim_{i\rightarrow \infty}E_{P}[\varphi(x_{j},B_{T}^{t,\theta^{ij}})|\mathcal{F}%
_{t}]=\lim_{i\rightarrow \infty}E_{P}[\varphi(x_{j},B_{T}^{t,\theta^{ij}}%
)]=\Lambda_{t,T}[x_{j}]=u_{t,T}(x_{j}).
\]
Setting $\theta^{i}=\sum_{j=1}^{N}\theta^{ij}I_{A_{j}}$, we have%
\begin{align*}
\Lambda_{t,T}[\zeta]  &  \geq E_{P}[\varphi(\zeta,B_{T}^{t,\theta^{i}%
})|\mathcal{F}_{t}]=E_{P}[\varphi(\sum_{j=1}^{N}I_{A_{j}}x_{j},B_{T}%
^{t,\sum_{j=1}^{N}I_{A_{j}}\theta^{ij}})|\mathcal{F}_{t}]\\
&  =\sum_{j=1}^{N}I_{A_{j}}E_{P}[\varphi(x_{j},B_{T}^{t,\theta^{ij}%
})|\mathcal{F}_{t}]\rightarrow \sum_{j=1}^{N}I_{A_{j}}u_{t,T}(x_{j}%
)=u_{t,T}(\zeta).
\end{align*}
On the other hand, for each given $\theta \in
\mathcal{A}_{t,T}^{\Theta}$, we have
\begin{align*}
E_{P}[\varphi(\zeta,B_{T}^{t,\theta})|\mathcal{F}_{t}]  &
=E_{P}[\varphi (\sum_{j=1}^{N}I_{A_{j}}x_{j},B_{T}^{t,\theta%
})|\mathcal{F}_{t}]\\
&  =\sum_{j=1}^{N}I_{A_{j}}E_{P}[\varphi(x_{j},B_{T}^{t,\theta%
})|\mathcal{F}_{t}]\\
&  \leq \sum_{j=1}^{N}I_{A_{j}}u_{t,T}(x_{j})=u_{t,T}(\zeta).
\end{align*}
We thus have $ess\sup_{\theta \in
\mathcal{A}_{t,T}^{\Theta}}E_{P}[\varphi
(\zeta,B_{T}^{t,\theta})|\mathcal{F}_{t}]\leq u_{t,T}(\zeta)$. The
proof is complete.
\end{proof}

The following result generalizes the well-known dynamical
programming principle:

\begin{theorem}
\label{Thm-GDDP}For each $\varphi \in C_{b.Lip}(\mathbb{R}^{n}\times \mathbb{R}^{2d}%
)$, $0\leq s\leq t\leq T$ and $\zeta \in L^{2}(\Omega,\mathcal{F}%
_{s},P;\mathbb{R}^{n})$ we have
\begin{equation}
ess\sup_{\theta \in \mathcal{A}_{s,T}^{\Theta}}E_{P}[\varphi(\zeta
,B_{t}^{s,\theta},B_{T}^{t,\theta})|\mathcal{F}_{s}]=ess\sup_{\theta
\in \mathcal{A}_{s,t}^{\Theta}}E_{P}[\psi(\zeta,B_{t}^{s,\theta})|\mathcal{F}%
_{s}], \label{eq5}%
\end{equation}
where $\psi \in C_{b.Lip}(\mathbb{R}^{n}\times \mathbb{R}^{d})$ is
given by
\[
\psi(x,y):=ess\sup_{\bar{\theta}\in \mathcal{A}_{t,T}^{\Theta}}E_{P}%
[\varphi(x,y,B_{T}^{t,\bar{\theta}})|\mathcal{F}_{t}]=\sup_{\bar{\theta}\in
\mathcal{A}_{t,T}^{\Theta}}E_{P}
[\varphi(x,y,B_{T}^{t,\bar{\theta}})].
\]

\end{theorem}

\begin{proof}
It is clear that
\[
ess\sup_{\theta \in \mathcal{A}_{s,T}^{\Theta}}E_{P}[\varphi(\zeta
,B_{t}^{s,\theta},B_{T}^{t,\theta})|\mathcal{F}_{s}]=ess\sup_{\theta
\in \mathcal{A}_{s,t}^{\Theta}}\left \{  ess\sup_{\bar{\theta}\in \mathcal{A}%
_{t,T}^{\Theta}}E_{P}[\varphi(\zeta,B_{t}^{s,\theta},B_{T}^{t,\bar{\theta}%
})|\mathcal{F}_{s}]\right \}  .
\]
It follows from (\ref{eq3}) and Lemma \ref{lem3} that
\begin{align*}
ess\sup_{\bar{\theta}\in
\mathcal{A}_{t,T}^{\Theta}}E_{P}[\varphi(\zeta
,B_{t}^{s,\theta},B_{T}^{t,\bar{\theta}})|\mathcal{F}_{s}]  &  =E_{P}%
[ess\sup_{\bar{\theta}\in
\mathcal{A}_{t,T}^{\Theta}}E_{P}[\varphi(\zeta
,B_{t}^{s,\theta},B_{T}^{t,\bar{\theta}})|\mathcal{F}_{t}]|\mathcal{F}_{s}]\\
&  =E_{P}[\psi(\zeta,B_{t}^{s,\theta})|\mathcal{F}_{s}],
\end{align*}
We thus have (\ref{eq5}).
\end{proof}

For each given $\varphi \in C_{b.Lip}(\mathbb{R}^{d})$ and $(t,x)\in
\lbrack 0,T]\times \mathbb{R}^{d}$, we set
\[
v(t,x):=\sup_{\theta \in \mathcal{A}_{t,T}^{\Theta}}E_{P}[\varphi(x+B_{T}%
^{t,\theta})].
\]
Since for each $h\in \lbrack0,T-t]$,
\begin{align*}
v(t,x)  &  =\sup_{\theta \in
\mathcal{A}_{t,T}^{\Theta}}E_{P}[\varphi
(x+B_{T}^{t,\theta})]\\
&  =\sup_{\theta \in \mathcal{A}_{t,T}^{\Theta}}E_{P}[\varphi(x+B_{t+h}%
^{t,\theta}+B_{T}^{t+h,\theta})]\\
&  =\sup_{\theta \in
\mathcal{A}_{t,t+h}^{\Theta}}E_{P}[v(t+h,x+B_{t+h}^{t,\theta })].
\end{align*}
This gives us the well-known dynamic programming principle:

\begin{proposition}
We have
\begin{equation}
v(t,x)=\sup_{\theta \in \mathcal{A}_{t,t+h}^{\Theta}}E_{P}[v(t+h,x+B_{t+h}%
^{t,\theta})]. \label{eq6}%
\end{equation}

\end{proposition}

\begin{lemma}
$v$ is bounded by $\sup|\varphi|$. It is a Lipschitz function in $x$
and $\frac{1}{2}$-h\"{o}lder function in $t$.
\end{lemma}

\begin{proof}
We only need to prove the regularity in $t$.
\[
\sup_{\theta \in
\mathcal{A}_{t,t+h}^{\Theta}}E_{P}[v(t+h,x+B_{t+h}^{t,\theta
})-v(t+h,x)]=v(t,x)-v(t+h,x).
\]
Since $v$ is a Lipschitz function in $x$, the absolute value of the
left hand
is bounded by%
\[
C\sup_{\theta \in
\mathcal{A}_{t,t+h}^{\Theta}}E_{P}[|B_{t+h}^{t,\theta}|]\leq
C_{1}h^{1/2}.
\]
The $\frac{1}{2}$-h\"{o}lder of $v$ in $t$ is obtained.
\end{proof}

\begin{theorem}
$v$ is a viscosity solution of the $G$-heat equation:%
\begin{align*}
\frac{\partial v}{\partial t}+G(D^{2}v)  &  =0,\  \  \text{on
}(t,x)\in
\lbrack0,T)\times \mathbb{R}^{d},\\
v(T,x)  &  =\varphi(x),
\end{align*}
where the function $G$ is given in (\ref{G}).
\end{theorem}

\begin{proof}
Let $\psi \in C_{b}^{2,3}((0,T)\times \mathbb{R}^{d})$ be such that
$\psi \geq v$ and, for a fixed $(t,x)\in(0,T)\times \mathbb{R}^{d}$,
$\psi(t,x)=v(t,x)$. From
the dynamic programming principle (\ref{eq6}) it follows that%
\begin{align*}
0  &  =\sup_{\theta \in \mathcal{A}_{t,t+h}^{\Theta}}E_{P}[v(t+h,x+B_{t+h}%
^{t,\theta})-v(t,x)]\\
&  \leq \sup_{\theta \in \mathcal{A}_{t,t+h}^{\Theta}}E_{P}[\psi(t+h,x+B_{t+h}%
^{t,\theta})-\psi(t,x)]\\
&  =\sup_{\theta \in \mathcal{A}_{t,t+h}^{\Theta}}E_{P}\left[  \int_{t}%
^{t+h}\left(  \frac{\partial \psi}{\partial
s}+\frac{1}{2}\text{tr}[\theta
_{s}\theta_{s}^{T}D^{2}\psi]\right)  (s,x+\int_{t}^{s}\theta_{r}dW_{r}%
)ds\right]  .
\end{align*}
Since $(\frac{\partial \psi}{\partial
s}+\frac{1}{2}$tr$[\theta_{s}\theta _{s}^{T}D^{2}\psi])(s,y)$ is
uniformly Lipschitz in $(s,y)$, we have for small $h>0$
\begin{align*}
&  E_{P}\mathbb{[}\left(  \frac{\partial \psi}{\partial s}+\frac{1}{2}%
\text{tr}[\theta_{s}\theta_{s}^{T}D^{2}\psi]\right)  (s,x+\int_{t}^{s}%
\theta_{r}dW_{r})]\\
&  \leq E_{P}[\left(  \frac{\partial \psi}{\partial s}+\frac{1}{2}%
\text{tr}[\theta_{s}\theta_{s}^{T}D^{2}\psi]\right)  (t,x)]+Ch^{1/2}.%
\end{align*}
Thus
\[
\sup_{\theta \in
\mathcal{A}_{t,t+h}^{\Theta}}E_{P}\int_{t}^{t+h}\left(
\frac{\partial \psi}{\partial s}+\frac{1}{2}\text{tr}[\theta_{s}\theta_{s}%
^{T}D^{2}\psi]\right)  (t,x)ds+Ch^{3/2}\geq0.
\]
Thus
\[
\left(  \frac{\partial \psi}{\partial s}+\frac{1}{2}\sup_{\gamma \in
\Theta }\text{tr}[\gamma \gamma^{T}D^{2}\psi]\right)
(t,x)h+Ch^{3/2}\geq0
\]
and then $[\frac{\partial \psi}{\partial
t}+G(D^{2}\psi)](t,x)\geq0$. By the definition, $v$ is a viscosity
subsolution. Similarly we can prove that it is also a supersolution.
\end{proof}

We observe that $u(t,x):=v(T-t,x)$, thus $u$ is a viscosity solution
of $\frac{\partial u}{\partial t}-G(D^{2}u)=0$, with Cauchy
condition $u(0,x)=\varphi(x)$.

From the uniqueness of the viscosity solution of $G$-heat equation
and Theorem \ref{Thm-GDDP}, we get immediately:

\begin{proposition}%
\label{prop50}
\begin{align*}
\mathbb{E}[\varphi(B_{t_{1}}^{0},B_{t_{2}}^{t_{1}},\cdots,B_{t_{n}}^{t_{n-1}%
})]&=\sup_{\theta \in \mathcal{A}_{0,T}^{\Theta}}E_{P}[\varphi(B_{t_{1}%
}^{0,\theta},B_{t_{2}}^{t_{1},\theta},\cdots,B_{t_{n}}^{t_{n-1},\theta})]\\
&=\sup_{\theta \in
\mathcal{A}_{0,T}^{\Theta}}E_{P_{\theta}}[\varphi(B_{t_{1}}^{0},B_{t_{2}}^{t_{1}},\cdots
,B_{t_{n}}^{t_{n-1}})],
\end{align*}
where $P_{\theta}$ is the law of the process $B_{t}^{0,\theta}=\int
_{0}^{t}\theta_{s}dW_{s}$, $t\geq0$, for $\theta \in
\mathcal{A}_{0,\infty }^{\Theta}$.
\end{proposition}

Now we prove that $\{P_{\theta},\  \theta \in
\mathcal{A}_{0,\infty}^{\Theta }\}$ is tight, this is important in
the following subsection.

\begin{proposition}
\label{prop51} The family of probability measures $\{P_{\theta}$:
$\theta \in \mathcal{A}_{0,\infty}^{\Theta}\}$ on
$C_{0}^{d}(\mathbb{R}^{+})$ is tight.
\end{proposition}

\begin{proof}
We apply It\^{o}'s formula to \ $(B_{t}^{s,\theta})_{t\geq s}$:%
\[
\left \vert B_{t}^{s,\theta}\right \vert ^{4}=\int_{s}^{t}4\left
\vert B_{r}^{s,\theta}\right \vert
^{2}B_{r}^{s,\theta}dB_{r}^{s,\theta}+2\int
_{s}^{t}\text{tr}[\theta_{r}\theta_{r}^{T}(I_{d}|B_{r}^{s,\theta}|^{2}%
+2B_{r}^{s,\theta}(B_{r}^{s,\theta})^{T})]dr.
\]
We thus have%
\begin{align*}
\mathbb{E}_{P_{\theta}}[|B_{t}-B_{s}|^{4}]  &  =2E\int_{s}^{t}\text{tr}%
[\theta_{r}\theta_{r}^{T}(I_{d}|B_{r}^{s,\theta}|^{2}+2B_{r}^{s,\theta}%
(B_{r}^{s,\theta})^{T})]dr\\
&  =2E\int_{s}^{t}[|\theta_{r}|^{2}|B_{r}^{s,\theta}|^{2}+2(\theta_{r}%
,B_{r}^{s,\theta})^{2}]dr\\
&  \leq C\int_{s}^{t}|B_{r}^{s,\theta}|^{2}dr\leq Cd\int_{s}^{t}(r-s)dr\\
&  =Cd\frac{(t-s)^{2}}{2}.
\end{align*}
We then can apply the well-known result of moment criterion for
tightness of Kolmogorov-Chentsov's type to conclude that
$\{P_{\theta}$: $\theta \in \mathcal{A}_{\Theta}\}$ is tight.
\end{proof}

%\begin{problem}
%$C_{lip}(\Omega)$ and $C_{b}(\Omega)$
%\end{problem}

\subsection{Capacity related to $G$-expectation}

We denote $\mathscr{P}_1=\{ P_{\theta}: \theta\in
\mathcal{A}_{0,\infty}^{\Theta} \}$ and
$\mathscr{P}=\overline{\mathscr{P}}_1$ the closure of
$\mathscr{P}_1$ under the topology of weak convergence. By
Proposition \ref{prop51}, $\mathscr{P}_1$ is tight and then
$\mathscr{P}$ is weakly compact. We set
\[
c(A):=\sup_{P\in\mathscr{P}}P(A),\  \  \ A\in\mathscr{B}(\Omega).
\]
For each $X\in L^0(\Omega)$ such that $E_{P}[X]$ exists for each
$P\in \mathscr{P}$, we set
\[
\hat{\mathbb{E}}[X]=\mathbb{E}^{\mathscr{P}}[X]=\sup_{P\in\mathscr{P}}E_{P}[X].
\]
Now we prove that
\begin{align*}
L^{1}_{G}(\Omega_T)&=\{ X\in L^0(\Omega_T): X\ \text{has a q.c.
version, }
\lim_{n\rightarrow\infty}\hat{\mathbb{E}}[|X|\mathbf{1}_{\{|X|>n\}}]=0
\},\\
L^{1}_{G}(\Omega)&=\{ X\in L^0(\Omega): X\ \text{has a q.c. version,
}
\lim_{n\rightarrow\infty}\hat{\mathbb{E}}[|X|\mathbf{1}_{\{|X|>n\}}]=0
\},\\
 \mathbb{E}[X]&=\hat{\mathbb{E}}[X], \  \  \ \forall X\in
L^{1}_{G}(\Omega),
\end{align*}
where q.c. denotes quasi-continuous for simplicity.

\  \  \
%We claim that the space $L_{ip}(\Omega)$ does not generate the
%topology of $\Omega$. Otherwise, taking open balls $G_n=B(0,n)$ for
%each $n\in \mathbb{N}$. For each $G_n$, if $X\in L_{ip}(\Omega)$
%such that $0\leq X\leq \mathbf{1}_{G_n}$, then it is easy to show
%that $X=0$, which yields $c(G_n)=0$ according to Section 3. Note
%that $c(\Omega)=\lim_{n}c(G_n)$, then we get $c(\Omega)=0$, which
%contradicts $c(\Omega)=1$. Thus $L_{ip}(\Omega)$ does not generate
%the topology of $\Omega$ and Theorem \ref{Thm8} cannot be directly
%applied. Fortunately we have following Lemma:
For proving this we need the following lemma.

\begin{lemma}
\label{Lemma53}Let $K$ be  a compact subset of $\Omega_T$ equipped
with the distance $\rho(\omega^1,\omega^2)=\max_{0\leq t\leq
T}|\omega^1_t-\omega^2_t|$. Then for each $\Phi \in
C_{b}(\Omega_T)$, there exists a sequence $\left \{ \Phi_{n}\right
\} _{n=1}^{\infty}\subset L_{ip}(\Omega_T)$ with $\left \Vert
\Phi_{n}\right \Vert _{\sup}\leq \left \Vert \Phi \right \Vert
_{\sup}$ such that $\Phi_{n}$ converges uniformly to $\Phi$ on $K$.
\end{lemma}

\begin{proof}This is just the consequence of the Stone-Weierstrass
theorem.
\end{proof}

\begin{theorem}
\label{Thm54}We have
\begin{align*}
L^{1}_{G}(\Omega_T)&=\{ X\in L^0(\Omega_T): X\ \text{has a q.c.
version, }
\lim_{n\rightarrow\infty}\hat{\mathbb{E}}[|X|\mathbf{1}_{\{|X|>n\}}]=0
\},\\
 \mathbb{E}[X]&=\hat{\mathbb{E}}[X], \  \  \ \forall X\in
L^{1}_{G}(\Omega_T). \end{align*}
\end{theorem}

\begin{proof}It follows from Proposition \ref{prop50} that
\[
\mathbb{E}[X]=\hat{\mathbb{E}}[X], \  \  \ \forall X\in
L_{ip}(\Omega_T).
\]
Thus $L^{1}_{G}(\Omega_T)$ can be seen as the completion of
$L_{ip}(\Omega_T)$ under the norm $\hat{\mathbb{E}}[|\cdot|]$. For
any fixed $\psi\in C_b(\Omega_T)$, since $\mathscr{P}$ is tight, we
have for each $n\in \mathbb{N}$, there exists a compact set
$K_n\subset\Omega_T$ such that $c(K_n^c)<\frac{1}{n}$. For this
$K_n$, by Lemma \ref{Lemma53}, there exists a $\varphi_n\in
L_{ip}(\Omega_T)$ such that
\[
\left \Vert \varphi_n \right \Vert _{\sup} \leq \left \Vert \psi
\right \Vert _{\sup} \quad \text{and} \quad \sup_{\omega\in
K_n}|\varphi_n(\omega)-\psi(\omega)|<\frac{1}{n}.
\]
Thus
\[
\hat{\mathbb{E}}[|\varphi_n-\psi|]\leq 2\left \Vert \psi \right
\Vert _{\sup}c(K_n^c)+\frac{1}{n}c(K_n)<(2\left \Vert \psi \right
\Vert _{\sup}+1)\frac{1}{n}\rightarrow0.
\]
It then follows that $C_b(\Omega_T)\subset L^{1}_{G}(\Omega_T)$, by
Theorem \ref{Thm8} we obtain the result. \end{proof}

\  \

\begin{remark}
The above results also hold for $L^{1}_{G}(\Omega)$, the proof is
similar.
\end{remark}

\  \

We also set
\[
\bar{c}(A):=\sup_{P\in\mathscr{P}_1}P(A),\  \  \
A\in\mathscr{B}(\Omega).
\]
It is easy to verify the following
\begin{enumerate}
  \item $\bar{c}(A)\leq c(A)$ for each $A\in\mathscr{B}(\Omega)$.
  \item $\bar{c}(O)= c(O)$ for each open set $O\subset\Omega$.
\end{enumerate}

Thus, a  function is $c$-quasi-continuous if and only if it is
$\bar{c}$-quasi-continuous, so we simply write quasi-continuous
function. For each $X\in L^0(\Omega)$ such that $E_{P}[X]$ exists
for each $P\in \mathscr{P}_1$, we set
\[
\bar{\mathbb{E}}[X]=\mathbb{E}^{\mathscr{P}_1}[X]=\sup_{P\in\mathscr{P}_1}E_{P}[X].
\]
It is easy to verify the following
\begin{enumerate}
  \item $\bar{\mathbb{E}}[X]\leq\hat{\mathbb{E}}[X]$ for each $X$
  which makes both expectation meaningful.
  \item $\bar{\mathbb{E}}[X]=\hat{\mathbb{E}}[X]$ for each
  bounded quasi-continuous function $X$.
\end{enumerate}
Similar to the proof of Theorem \ref{Thm54}, we get the following
theorem.
\begin{theorem}
We have
\begin{align*}
L^{1}_{G}(\Omega_T)&=\{ X\in L^0(\Omega_T): X\ \text{has a q.c.
version, }
\lim_{n\rightarrow\infty}\bar{\mathbb{E}}[|X|\mathbf{1}_{\{|X|>n\}}]=0
\},\\
L^{1}_{G}(\Omega)&=\{ X\in L^0(\Omega): X\ \text{has a q.c. version,
}
\lim_{n\rightarrow\infty}\bar{\mathbb{E}}[|X|\mathbf{1}_{\{|X|>n\}}]=0
\},\\
\mathbb{E}[X]&=\bar{\mathbb{E}}[X], \  \  \ \forall X\in
L^{1}_{G}(\Omega).
\end{align*}
\end{theorem}

\begin{remark}
Theorem \ref{Thm15} holds for $\hat{\mathbb{E}}[\cdot]$ under the
capacity $c(\cdot)$. But it does not necessarily hold for
$\bar{\mathbb{E}}[\cdot]$ under the capacity $\bar{c}(\cdot)$.
\end{remark}

\  \

\section*{Acknowledgement} This works was initiated while L. Denis was visiting the
Shandong University in Jinan. He wishes to thank Pr. Peng and all
the members of the School of Mathematics in Jinan for their kind
hospitality. His work is supported by the chair "risque de cr\'edit", F\'ed\'eration bancaire Fran\c{c}aise.\\
S. Peng thanks the partial support from The National
Basic Research Program of China (973 Program) grant No. 2007CB814900
(Financial Risk).

\end{document}